\documentclass[12pt,reqno]{amsart}
\usepackage{amsmath,amssymb,amsfonts,amscd,latexsym,amsthm,mathrsfs,comment}
\usepackage[usenames]{color}
\usepackage[unicode]{hyperref}
\theoremstyle{remark}

\begin{document}

\title[A supercongruence fantasy on Fibonacci and Lucas]{A supercongruence fantasy on Fibonacci and Lucas\\(after Guillera)}

\author{Wadim Zudilin}
\address{Department of Mathematics, IMAPP, Radboud University, PO Box 9010, 6500~GL Nij\-me\-gen, Netherlands}
\urladdr{https://www.math.ru.nl/~wzudilin/}

\date{5 February 2026}

\dedicatory{To Jes\'us Guillera for numerous mathematical inspirations\,---\,your work deserves a prize!}


\maketitle

\begin{abstract}
Motivated by observations of Guillera we generalise the so-called Ramanujan-type supercongruences to a further level in which the sequences of Fibonacci, Lucas, Ap\'ery numbers and their friends all receive a natural appearance.
\end{abstract}

Among the seventeen series for $1/\pi$ recorded by Ramanujan \cite{Ra14} in 1914, there was just one `irrational' entry
\begin{equation}
\sum_{n=0}^\infty{\binom{2n}n}^3\big((30+45\sqrt5)n+(-1+5\sqrt5)\big)\frac{(-1/\varphi)^{8n}}{2^{12n}}
=\frac{32}\pi,
\label{eq1}
\end{equation}
where $\varphi=(1+\sqrt5)/2$ denotes the golden ratio; note that $-1/\varphi=(1-\sqrt5)/2$ is the algebraic conjugate of~$\varphi$.
A century later, in 2020, Guillera \cite{Gu20} observed that combining the formula with its `conjugate companion'
\begin{equation}
\sum_{n=0}^\infty{\binom{2n}n}^3\big((30-45\sqrt5)n+(-1-5\sqrt5)\big)\frac{\varphi^{8n}}{2^{12n}}
=-\frac{96}\pi,
\label{eq2}
\end{equation}
(whose convergence is not as fast as that of \eqref{eq1}) leads to a pair of two `rational' Ramanujan-like series
\begin{align}
\sum_{n=0}^\infty{\binom{2n}n}^3\frac{(-30F_{8n}+42L_{8n})n+(F_{8n}+5L_{8n})}{2^{12n}}
&=\frac{128}{\pi\sqrt5},
\label{eq3}
\\
\sum_{n=0}^\infty{\binom{2n}n}^3\frac{(210F_{8n}-30L_{8n})n+(25F_{8n}+L_{8n})}{2^{12n}}
&=\frac{64}{\pi};
\label{eq4}
\end{align}
here $F_n$ and $L_n$ are the Fibonacci and Lucas numbers, respectively, generated by the expansions
\[
\sum_{n=0}^\infty F_nt^n=\frac t{1-t-t^2}=t+t^2+\dotsb
\quad\text{and}\quad
\sum_{n=0}^\infty L_nt^n=\frac{2-t}{1-t-t^2}=2+t+3t^2+\dotsb.
\]
Behind the 1-page scenes in \cite{Gu20}, there are congruences
\begin{equation}
\begin{gathered}
S_1(p)\overset?\equiv 10\bigg(\frac{-5}p\bigg)p\pmod{p^3}
\quad\text{and}\quad
S_2(p)\overset?\equiv 2\bigg(\frac{-1}p\bigg)p\pmod{p^3},
\\
\text{where $p>2$ is a prime},
\end{gathered}
\label{eq5}
\end{equation}
for the truncated sums
\begin{align*}
S_1(N)&=\sum_{n=0}^{N-1}{\binom{2n}n}^3\frac{(-30F_{8n}+42L_{8n})n+(F_{8n}+5L_{8n})}{2^{12n}},
\\
S_2(N)&=\sum_{n=0}^{N-1}{\binom{2n}n}^3\frac{(210F_{8n}-30L_{8n})n+(25F_{8n}+L_{8n})}{2^{12n}}
\end{align*}
and the Legendre symbol $\big(\frac{m}{p}\big)$.
This type of Ramanujan supercongruences is discussed by Guillera in \cite{Gu12}, and we now expect, more generally, the validity of their generalisations \cite{GZ21,Sw15}
\begin{equation}
S_1(p^s)\overset?\equiv\bigg(\frac{-5}p\bigg)pS_1(p^{s-1})\pmod{p^{3s}}
\quad\text{and}\quad
S_2(p^s)\overset?\equiv\bigg(\frac{-1}p\bigg)pS_2(p^{s-1})\pmod{p^{3s}},
\label{eq6}
\end{equation}
as well as
\begin{equation}
\begin{gathered}
S_1\big((p^s+1)/2\big)\overset?\equiv\bigg(\frac{-5}p\bigg)pS_1\big((p^{s-1}+1)/2\big)\pmod{p^{3s}}
\quad\text{and}\\
S_2\big((p^s+1)/2\big)\overset?\equiv\bigg(\frac{-1}p\bigg)pS_2\big((p^{s-1}+1)/2\big)\pmod{p^{3s}},
\end{gathered}
\label{eq7}
\end{equation}
where $p>2$ is a prime and $s=1,2,\dots$\,.

In a very much the same spirit, one can `rationalise' the Ramanujan--Guillera-type formula
\begin{equation}
\sum_{n=0}^\infty{\binom{2n}n}^4\binom{3n}{n}\big(20(61-27\sqrt5)n^2+9(101-45\sqrt5)n+3(56-25\sqrt5)\big)\frac{(1/\varphi)^{15n}}{2^{6n}}
\overset?=\frac{3}{\pi^2},
\label{eq8}
\end{equation}
discovered by Almkvist and Guillera \cite{AG12}, and associate with it a formal conjugate companion (which now diverges!) and the related truncated sums
{\small
\begin{align*}
S_3(N)&=\sum_{n=0}^{N-1}{\binom{2n}n}^4\binom{3n}{n}(-1)^n\frac{(2440F_{15n} + 1080L_{15n})n^2 + (1818F_{15n} + 810L_{15n})n + (336F_{15n} + 150L_{15n})}{2^{6n}},
\\
S_4(N)&=\sum_{n=0}^{N-1}{\binom{2n}n}^4\binom{3n}{n}(-1)^n\frac{(5400F_{15n} + 2440L_{15n})n^2 + (4050F_{15n} + 1818L_{15n})n + (750F_{15n} + 336L_{15n})}{2^{6n}}.
\end{align*}%
}
As before, we numerically observe the supercongruences
\begin{equation}
S_3(p^s)\overset?\equiv\bigg(\frac{5}p\bigg)p^2S_3(p^{s-1})\pmod{p^{5s}}
\quad\text{and}\quad
S_4(p^s)\overset?\equiv p^2S_4(p^{s-1})\pmod{p^{5s}},
\label{eq9}
\end{equation}
while the analogous ones for shorter truncations are slightly weaker:
\begin{equation}
\begin{gathered}
S_3\big((p^s+1)/2\big)\overset?\equiv\bigg(\frac{5}p\bigg)p^2S_3\big((p^{s-1}+1)/2\big)\pmod{p^{4s}}
\quad\text{and}\\
S_4\big((p^s+1)/2\big)\overset?\equiv p^2S_4\big((p^{s-1}+1)/2\big)\pmod{p^{4s}},
\end{gathered}
\label{eq10}
\end{equation}
where $p>2$ is a prime and $s=1,2,\dots$\,.

As a reasonable simplification for numerical check we include the generating functions for the rarefied Fibonacci and Lucas sequences involved in the hypergeometric sums $S_1(N)$, $S_2(N)$, $S_3(N)$ and $S_4(N)$:
\[
\sum_{n=0}^\infty F_{8n}t^n=\frac{21t}{1-47t+t^2}
\quad\text{and}\quad
\sum_{n=0}^\infty L_{8n}t^n=\frac{2-47t}{1-47t+t^2},
\]
together with
\[
\sum_{n=0}^\infty F_{15n}t^n=\frac{610t}{1-1364t-t^2}
\quad\text{and}\quad
\sum_{n=0}^\infty L_{15n}t^n=\frac{2-1364t}{1-1364t-t^2}.
\]
Note that such Fibonacci-like sequences satisfy recurrence relations (with constant rational coefficients) of order higher than~1 in contrast with the sequence of pure powers $\{z^n\}_{n\ge0}$ for a fixed rational $z\ne0$: the latter one solves such a recursion of order~1.
This indicates a principal difference between the classical and Fibonacci-like versions of Ramanujan-type formulas and supercongruences.
Notice however that the hypergeometric input
\[
\frac1{2^{12n}}{\binom{2n}n}^3 \quad\text{and}\quad 
\frac{(-1)^n}{2^{6n}}{\binom{2n}n}^4\binom{3n}{n}
\]
in the formulas above also comes from sequences satisfying linear recursions, this time with coefficients that are rational functions of~$n$, and of order~1.
It is already known (see, e.g., \cite{AG12,CZ10}) that such binomial-product expressions can be replaced with other arithmetic sequences of `Ap\'ery type' that are solutions of similar recursions but of order higher than~1; this kind of formulas for $1/\pi$ and its powers are dubbed in the literature as Ramanujan--Sato-type formulas.

One can generate a further mix of both kinds of these formulas, for example, from `rationalisation' of the Ramanujan--Sato-type formula \cite[eq.~(6.2)]{CZ10}
\[
\sum_{n=0}^\infty A_n(8n+4-\sqrt6)\rho^{2n+1}=\frac1{\pi\sqrt2},
\]
where $\rho=(\sqrt3-\sqrt2)^2=5-2\sqrt6$ is a root of the polynomial $1-10t+t^2$ and
\[
A_n=\sum_{k=0}^n{\binom{n}{k}}^2{\binom{n+k}{k}}^2
\]
are the Ap\'ery numbers. For this we introduce the Fibonacci-like numbers
\[
U_n=\frac{(5+2\sqrt6)^{2n} - (5-2\sqrt6)^{2n}}{40\sqrt6}
\quad\text{and}\quad
V_n=\frac{(5+2\sqrt6)^{2n} + (5-2\sqrt6)^{2n}}{2}
\]
with their generating functions
\[
\sum_{n=0}^\infty U_nt^n=\frac t{1-98t+t^2},
\quad
\sum_{n=0}^\infty V_nt^n=\frac{1-49t}{1-98t+t^2},
\]
and the related truncated companion sums
\begin{align*}
S_5(N)&=\sum_{n=0}^{N-1}A_n\big((800U_n + 16V_n)n + (640U_n + 13V_n)\big),
\\
S_6(N)&=\sum_{n=0}^{N-1}A_n\big((1920U_n + 40V_n)n + (1560U_n + 32V_n)\big).
\end{align*}
Then
\begin{equation}
S_5(p^s)\overset?\equiv\bigg(\frac{-3}p\bigg)pS_5(p^{s-1})\pmod{p^{3s}}
\quad\text{and}\quad
S_6(p^s)\overset?\equiv\bigg(\frac{-2}p\bigg)pS_6(p^{s-1})\pmod{p^{3s}}
\label{eq11}
\end{equation}
for all primes $p>2$ and $s=1,2,\dots$, with the sole exception of the case $S_6(3^1)$.

Such congruences all follow the pattern initially observed in \cite{Zu09}, \cite{Gu12} and \cite{Sw15}, and there was a considerable amount of work in the latest years to cover particular instances of cases of `pure hypergeometric and rational' type\,---\,nothing really for Fibonacci and friends to benefit from.
It would be really nice to have a unifying conceptual framework for understanding this quite general arithmetic phenomenon.

\end{document}